\documentclass[conference]{IEEEtran}
\IEEEoverridecommandlockouts
\usepackage{soulutf8}
\usepackage{cite}
\usepackage{amsmath,amssymb,amsfonts,bm}
\usepackage{graphicx}
\usepackage{xcolor}
\def\BibTeX{{\rm B\kern-.05em{\sc i\kern-.025em b}\kern-.08em
    T\kern-.1667em\lower.7ex\hbox{E}\kern-.125emX}}

\newcommand{\image}{\boldsymbol{x}}
\newcommand{\gt}{\bar{\image}}
\newcommand{\sm}{\boldsymbol{A}}
\newcommand{\nfd}{\bar{\boldsymbol{b}}}
\newcommand{\noise}{\boldsymbol{e}}
\newcommand{\rhs}{\boldsymbol{b}}
\newcommand{\DD}{\boldsymbol{D}}
\newcommand{\MM}{\boldsymbol{M}}
\newcommand{\II}{\boldsymbol{I}}

\newcommand{\EE}{\mathcal{E}}
\newcommand{\Var}{\mathrm{Cov}}

\newcommand{\residual}{\boldsymbol{\varrho}}
\newcommand{\relaxp}{\omega}

\newcommand{\row}{\boldsymbol{r}}

\newcommand{\vv}{\boldsymbol{v}}
\newcommand{\ww}{\boldsymbol{w}}
\newcommand{\boldzero}{\bm{0}}
\newcommand{\vz}{\boldsymbol{z}}
\newcommand{\ncp}{\boldsymbol{c}}

\begin{document}

\title{Stopping Rules for Algebraic Iterative
 Reconstruction Methods in Computed Tomography\thanks{This work was
 partially funded by a
 Villum Investigator grant (no. 25893) from The Villum Foundation.}
}

\author{\IEEEauthorblockN{Per Christian Hansen, Jakob Sauer J{\o}rgensen,
 Peter Winkel Rasmussen}
\IEEEauthorblockA{\textit{Department of Applied Mathematics and Computer Science} \\
\textit{Technical University of Denmark},
Kgs. Lyngby, Denmark \\
ORCID 0000-0002-7333-7216, 0000-0001-9114-754X, 0000-0002-0823-0316 \\
Email \texttt{\{pcha,jakj,pwra\}@dtu.dk} }
}

\maketitle

\begin{abstract}
Algebraic models for the reconstruction problem in X-ray computed tomography (CT)
provide a flexible framework that applies to many measurement geometries.
For large-scale problems we need to use iterative solvers, and we need
stopping rules for these methods that terminate the iterations when we
have computed a satisfactory reconstruction that balances the
reconstruction error and the influence of noise from the measurements.
Many such stopping rules are developed in the inverse problems communities,
but they have not attained much attention in the CT world.
The goal of this paper is to describe and illustrate four stopping rules
that are relevant for CT reconstructions.
\end{abstract}

\begin{IEEEkeywords}
tomographic reconstruction, iterative methods, stopping rules, semi-convergence
\end{IEEEkeywords}

\section{Introduction}

This paper considers large-scale methods for computed tomographic (CT)
based on a discretization of the problem followed by solving the
system of linear equations by means of an iterative solver.
These methods are quite generic in the sense
that they do not assume any specific scanning geometry, and they tend to produce
good reconstructions with few artifacts
in the case of limited-data and/or limited-angle problems.

In CT, a forward projection
maps the object to the data in the form of
projections of the object onto the detector planes for various scan positions.
In the case of parallel-beam CT the forward projection is known as the
Radon transform~\cite{Natterer}.

In practise, data consists of noisy measurements of the attenuation of the X-rays
that pass through the object.
The discretization of the reconstruction problem takes the form
  \begin{equation}
    \label{eq:Ax}
    \sm\, \image \approx \rhs \ , \qquad \rhs = \sm\,\gt + \noise \ ,
  \end{equation}
where the ``system matrix'' $\sm \in \mathbb{R}^{m\times n}$ is a discretization
of the forward projector,
$\rhs\in\mathbb{R}^m$ is a vector with the measured data,
and $\image\in\mathbb{R}^n$ is a vector that holds the pixels of the
reconstructed image of the object's interior.
Moreover, $\gt$ represents the exact object and $\noise$
represents the measurement noise.
A number of discretization schemes are available for computing the matrix $\sm$,
see, e.g., \cite{SIAMbook,KSSHN}.

There are no restrictions on the dimensions $m$ and $n$ of the matrix~$\sm$; both
over-determined and under-determined systems are common, depending on the
measurement setup.
The matrix $\sm^T$ represents the so-called back projector which
maps the data back onto the solution domain~\cite{Natterer};
it plays a central role in filtered back projection
and similar methods.

In large-scale CT  problems, the matrix $\sm$ -- in spite of
the fact that it is sparse -- is often too large to store explicitly.
For this reason we must use iterative solvers that only access
the matrix via functions that compute the multiplications with $\sm$
and $\sm^T$ in a matrix-free fashion, often using GPUs or other
hardware accelerators.
In CT these iterative solvers are collectively referred to as
\emph{algebraic iterative reconstruction methods} which includes
well-known methods such as ART \cite{GoBH70} and SIRT
(also known as SART) \cite{SIRT}.

\begin{figure*}
  \centering
  \includegraphics[width=0.68\textwidth]{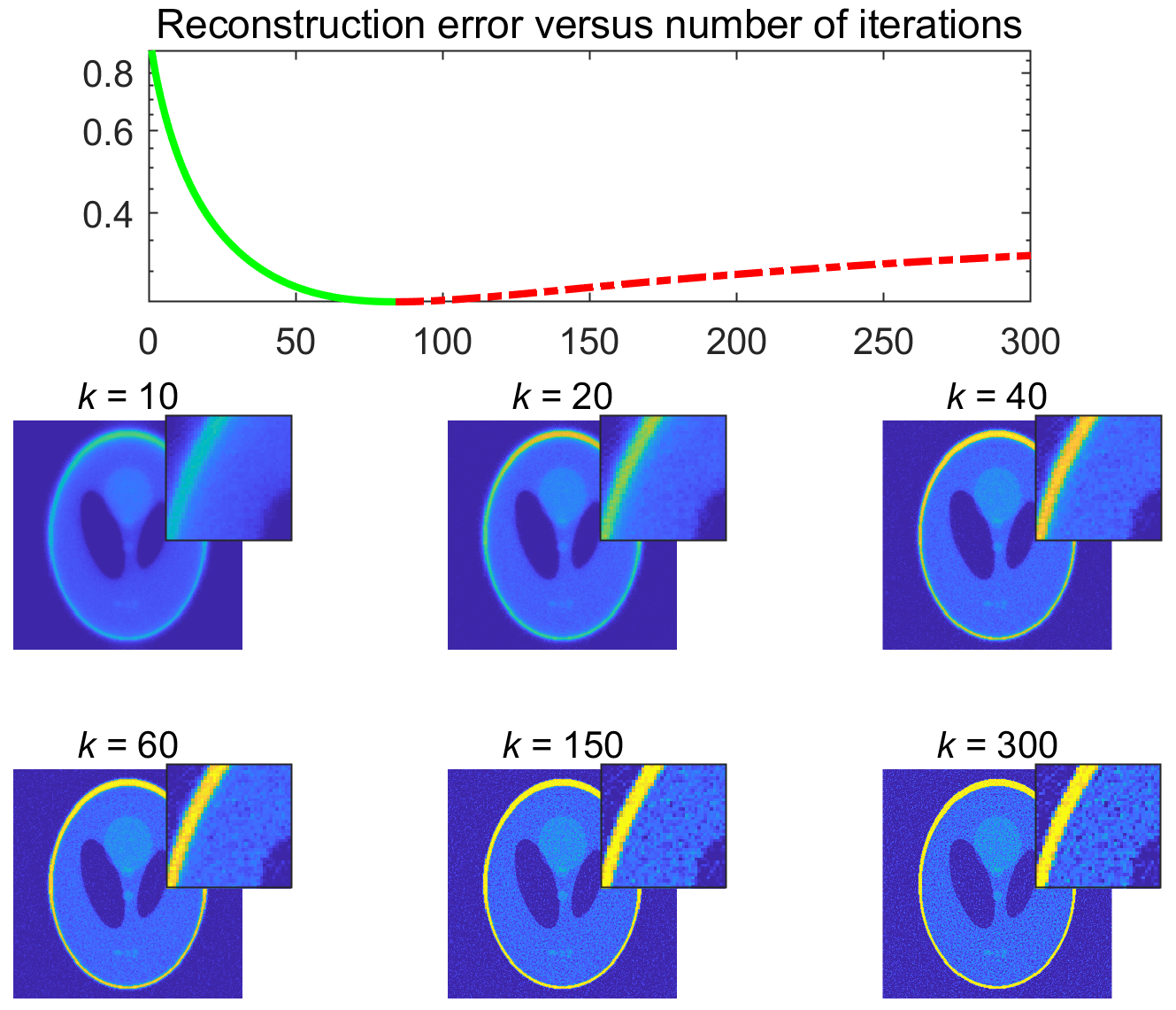}
  \caption{\label{fig:IterationImagesKaczmarz}Illustration of semi-convergence
  for Landweber's method applied to a small noisy test problem.
  Top:\ the error history, i.e., the reconstruction error
  $\| \image^{(k)} - \gt \|_2$ as a function of the number of iterations~$k$.
  The solid green line shows the part where $\image^{(k)}$ approaches $\gt$
  while the red dash-dotted line shows the part where $\image^{(k)}$
  becomes dominated by the noise.
  Bottom:\ selected iterations with inserts that zoom in on a small region;
  the increasing amount of noise is clearly visible.}
\end{figure*}

Common for all these methods is that they, from an initial vector
$\image^{(0)}$ (often the zero vector) produce a sequence of iteration vectors
$\image^{(k)}$, $k=1,2,\ldots$ which, in the ideal situation,
converge to the ground truth~$\bar{x}$.
In practise, however, when noise is present in the measured data
we experience a phenomenon called \emph{semi-convergence:}
\begin{itemize}
\item
During the initial iterations, the iteration vector $\image^{(k)}$
approaches the desired but un-obtainable ground truth~$\gt$. \\[-2.4mm]
\item
During later iterations, $\image^{(k)}$ converges to the undesired noisy
solution associated with the particular iterative method
(e.g., $\sm^{-1}\rhs$ 
if the system matrix is invertible).
\end{itemize}
This is illustrated in Fig.~\ref{fig:IterationImagesKaczmarz} which
shows the error history, i.e., the reconstruction error
$\| \image^{(k)} - \gt \|_2$ versus the number of iterations $k$,
together with selected iterates $\image^{(k)}$ shown as images.
The error history has the characteristic form associated with semi-convergence.

If we can stop the iterations just when the convergence behavior changes
from the former to the latter, then we achieve an approximation to
$\gt$ that is not too perturbed by the noise in the data.
\emph{This paper describes four such methods based on certain statistical
properties of the noise.}

Sections \ref{sec:four} and \ref{subsec:trace} describe four
state-of-the-art stopping rules as well as two methods to efficiently
estimate a trace-term that is needed on some of these rules; all
numerical experiments in these sections were performed by means of
the AIR Tools II software package~\cite{AIR}.
In Section \ref{sec:large} we illustrate these techniques with a
large-scale example.


\section{Four Stopping Rules}
\label{sec:four}

Our stopping rules apply to methods of the general form
  \begin{equation}
  \label{eq:SIRT}
    \image^{(k+1)} = \image^{(k)} + \DD \,\sm^T \MM
    \bigl( \rhs - \sm \, \image^{(k)} \bigr) \ ,
  \end{equation}
where different choices of the diagonal matrices $\DD$ and $\MM$ lead to
different methods -- see, e.g., \cite{Elfving} for an overview.
To simplify the presentation, we focus on the simple case where
$\DD$ and $\MM$ are identity matrices, in which case we obtain
\emph{Landweber's method} (which is equivalent to the steepest descent
method applied to the least squares problem):

\medskip\noindent
\framebox[0.48\textwidth][l]{\vbox{
  \begin{tabbing}
  xxx \= xxx \= xxx \kill
  \> \underline{Landweber's method} \\[2mm]
  \> $\image^{(0)}$ = initial vector \\
  \> for $k=0,1,2,\ldots$ \\[1mm]
  \> \> $\image^{(k+1)} = \image^{(k)} + \relaxp\,
     \sm^T (\rhs - \sm\,\image^{(k)} ) $ \\
  \> end
  \end{tabbing}
}}

\medskip

We will frequently refer to the \emph{residual}
for the $k$th iterate, defined as the vector
  \begin{equation}
  \label{eq:residual}
    \residual^{(k)} = \rhs - \sm\,\image^{(k)}\ , \qquad k=1,2,3,\ldots
  \end{equation}
Moreover, after a bit of algebraic manipulations it turns out that we
can write the $k$th iteration vector as
  \begin{equation}
  \label{eq:Aksharp}
    \image^{(k)} =
    \sum_{j=0}^{k-1} ( \II - \relaxp\, \sm^T \sm )^j \relaxp\,\sm^T\rhs
    = \sm_k^{\#} \rhs \ ,
  \end{equation}
which defines the matrix $\sm_k^{\#}$ such that we
can write the $k$th iterate as $\image^{(k)} = \sm_k^{\#} \rhs$.
This is convenient as a theoretical tool, but $\sm_k^{\#}$
is never computed explicitly.

To set the stage, we need to introduce
a small amount of statistical framework and notation.
We will often need the exact \emph{noise-free data} that corresponds to
the ground truth image, and we write
  \begin{equation}
  \label{eq:nfd}
    \nfd = \sm\,\gt \ .
  \end{equation}
We can then write the measured data as $\rhs = \nfd+\noise$.
The elements of the noise vector $\noise$ are random variables, i.e.,
their values depend on a set of well-defined random events.
The vector of expected values $\EE(\noise)$ and the covariance matrix
$\Var(\noise)$ are defined as
  \begin{align}
    \EE(\noise) &= \begin{pmatrix} \EE(e_1) \\ \EE(e_2) \\ \vdots \end{pmatrix} , \\[1mm]
    \Var(\noise) &= \EE \!\left( \bigl( \noise - \EE(\noise) \bigr)\,
      \bigl( \noise - \EE(\noise) \bigr)^T \right) .
  \end{align}
To simplify our discussion and make the ideas clearer, throughout this section
we will restrict our analysis to \emph{white Gaussian noise} with zero mean:
  \begin{equation}
  \label{eq:whiteGaussian}
    \EE(\noise) = \boldsymbol{0} \ , \quad \Var(\noise) = \eta^2\II \  ,
    \quad \EE(\|\noise\|_2^2) = m\,\eta^2 \ ,
  \end{equation}
where $\eta$ is the standard deviation of the noise and $m$ is the
number of elements in~$\noise$.
Noise in tomographic problems is rarely strictly Gaussian,
but sometimes this is a reasonable assumption.

\subsection{Fitting to the Noise Level}
\label{sec:fittonoiselevel}

Our description of this stopping rule is based on \cite[\S 11.2.3]{SIAMbook}.
A simple idea is to choose the number of iterations $k$
such that the residual $\residual^{(k)}$
is ``of the same size'' as the noise vector~$\noise$.
Specifically, such that $\| \residual^{(k)} \|_2$ approximates the expected value
$\EE(\|\noise\|_2)$ of the latter:
  \begin{equation}
  \label{eq:discrepancyprinciple}
    \| \residual^{(k)} \|_2 \approx \eta\,\sqrt{m} \ .
  \end{equation}
In the literature this is referred to as the
\emph{discrepancy principle} \cite{EnHN00}.
Since $\|\residual^{(k)}\|_2$ takes discrete values for
$k=1,2,3,\ldots$ we cannot expect to find a $k$ such
that the above holds with equality.

It is common to include a constant $\tau$ slightly larger
than 1, say, $\tau = 1.02$,
such that the above condition takes the form
$\| \residual^{(k)} \|_2 \leq \tau\,\eta\,\sqrt{m}$.
This constant can be useful as a ``safety factor''
when we have only a rough estimate of the noise.

If we replace $\image^{(k)}$ with the ground truth $\gt$ then the residual is
$\rhs - \sm\,\gt = \noise$ and the residual norm obviously equals $\|\noise\|_2$.
However, this is not a sound statistical argument that the norm of the residual
$\residual^{(k)}$ in Eq.\ (\ref{eq:residual}) should be equal to~$\|\noise\|_2$
for the optimal iterate $\image^{(k)}$.

Here we present an alternative that is based on statistical principles.
To motivate this stopping rule, we split the residual vectors as follows:
  \begin{align*}
    \residual^{(k)} & = \rhs - \sm\,\image^{(k)} = \rhs - \sm\, \sm_k^{\#} \rhs \\[1mm]
    &= (\II - \sm\, \sm_k^{\#}) \, \nfd \ + (\II - \sm\, \sm_k^{\#}) \, \noise \ .
  \end{align*}
The heuristic insight is then as follows:
\begin{itemize}
\item
When $k$ is too small then $\sm\,\image^{(k)}$ is not a good approximation the
exact data $\nfd$.  Hence, the residual $\residual^{(k)}$ is dominated by
$(\II - \sm\, \sm_k^{\#}) \, \nfd$
and $\| (\II - \sm\, \sm_k^{\#}) \, \nfd \|_2$ is
larger than $\| (\II - \sm\, \sm_k^{\#}) \, \noise \|_2$. \\[-2.4mm]
\item
When $k$ is ``just about right'' then $\sm\,\image^{(k)}$ approximates
$\nfd$ as well as possible; the norm
$\| (\II - \sm\, \sm_k^{\#}) \, \nfd \|_2$ has now become smaller
and it is of the same
size as the norm $\| (\II - \sm\, \sm_k^{\#}) \, \noise \|_2$. \\[-2.4mm]
\item
When $k$ is too large then the residual $\residual^{(k)}$ is dominated by the
noise component $(\II - \sm\, \sm_k^{\#}) \, \noise$,
and therefore $\| (\II - \sm\, \sm_k^{\#}) \, \noise \|_2$ dominates the residual norm.
\end{itemize}
According to these observations we should therefore choose $k$ such that
$\| (\II - \sm\, \sm_k^{\#}) \, \nfd \|_2 \approx \| (\II - \sm\, \sm_k^{\#}) \, \noise \|_2$.
Unfortunately both these are unknown.

The above heuristic reasoning has been formalized
in \cite{HaTi87}, \cite{TBKT91} and \cite{Turchin67}, and we
will summarize the main results as they apply here.
The key points are that $\| (\II - \sm\, \sm_k^{\#}) \, \rhs \|_2$
is an approximation
to the prediction error $\| (\II - \sm\, \sm_k^{\#}) \, \rhs \|_2$  and that
  \[
    \EE(\| (\II - \sm\, \sm_k^{\#}) \, \noise \|_2^2) = \eta^2 \, ( m - t_k)
  \]
in which
  \begin{equation}
  \label{eq:tk}
     t_k = \mathrm{trace}(\sm\sm_k^{\#}) \ .
  \end{equation}
Hence, at the optimal $k$ we have
  \begin{equation}
  \label{eq:EEresidual}
    \EE(\| \residual^{(k)} \|_2^2 ) \approx \eta^2 \, ( m - t_k) \ .
  \end{equation}
  Here, $k$ is ``optimal'' in the sense that it is the largest iteration number for which
we cannot reject $\image^{(k)}$ -- computed from the noisy data $\rhs$ -- as a possible
solution to the noise-free system, cf.\ \cite[p.~93]{Turchin67}.

The real number $m-t_k$ is sometimes referred to as the ef\-fective (or
equivalent) degrees of freedom \cite{Wahba90} in the residual.
An exact computation of $t_k$ is cumbersome for most
methods, but it can be
approximated quite efficiently as described in \S\ref{subsec:trace}.
We have thus arrived at the following
``fit-to-noise-level'' (FTNL) stopping rule where, again,
we include the ``safety factor'' $\tau$:

\medskip\noindent
\framebox[0.48\textwidth][l]{\vbox{
  \begin{tabbing}
  xxx \= xxx \= xxx \kill
  \> \underline{Stop rule:\ FTNL\,} \\[2mm]
  \>  Stop at the smallest $k$ \\[2mm]
  \> \> for which $\|\residual^{(k)} \|_2 \leq \tau \, \eta\,\sqrt{m-t_k} \ .$
  \end{tabbing}
}}

\medskip

\begin{figure}
  \centering
  \includegraphics[width=0.48\textwidth]{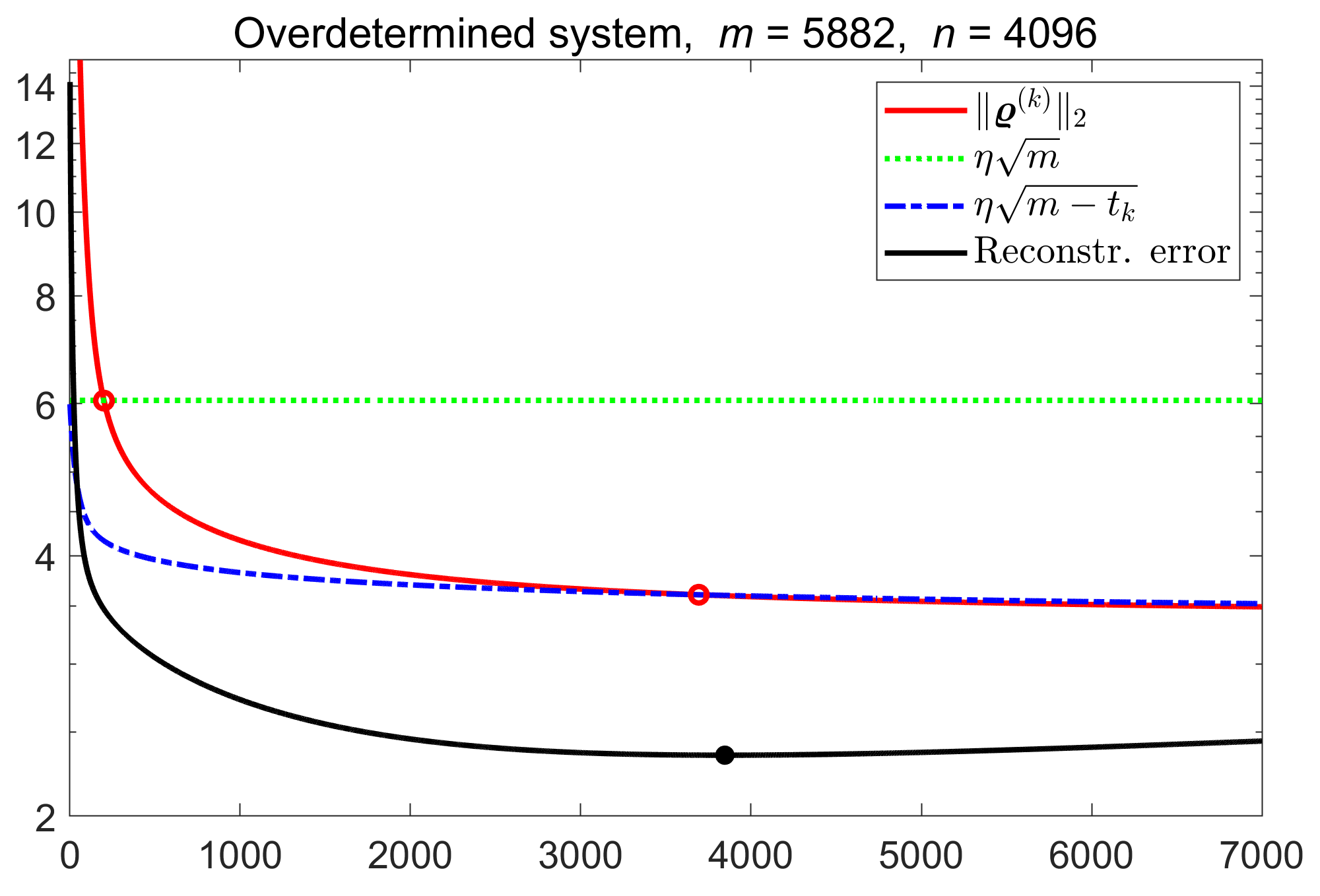}
  \includegraphics[width=0.48\textwidth]{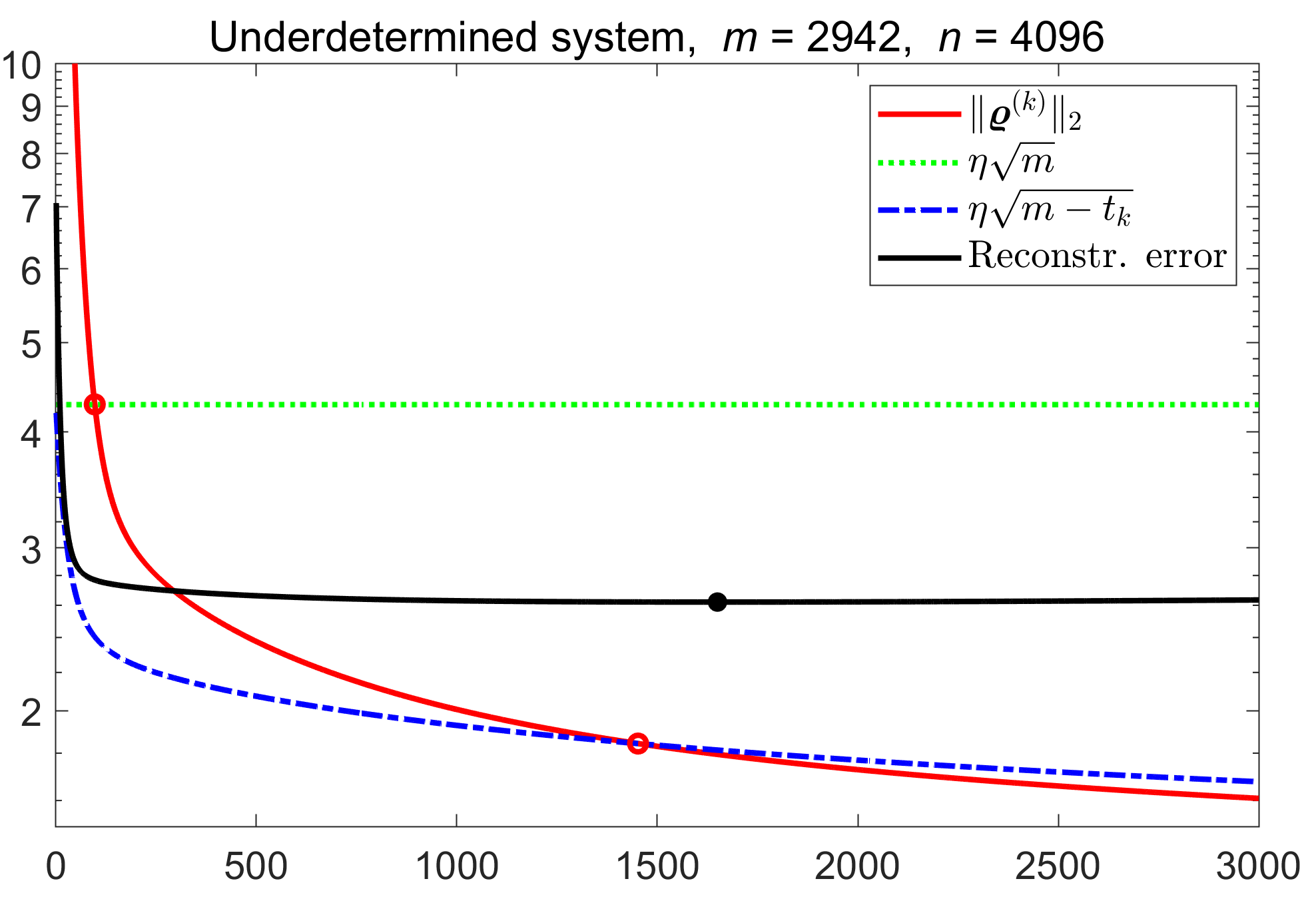}
  \caption{\label{fig:FitToNoise}Illustration of the FTNL
  stopping rule for Landweber's method,
  with two parallel-beam tomographic problems.
  The smallest reconstruction error is marked with the black dot,
  and the residual norms that satisfies the stopping rules are marked with
  red circles.  The FTNL rule works well, while
  stopping at that $k$ for which
  $\| \residual^{(k)} \|_2 \approx \eta\, \sqrt{m}$ terminates the iterations
  much too early.}
\end{figure}

\textbf{Example 1.}
We illustrate the FTNL ``fit-to-noise-level'' stopping rule with two small
parallel-beam CT problems with image size $64 \times 64$ and
$91$ detector pixels.
The projection angles are, respectively,
$3^{\circ},6^{\circ},9^{\circ},\ldots,180^{\circ}$ (giving an
over-determined system)
and $8^{\circ},16^{\circ},24^{\circ},\ldots,180^{\circ}$
(giving an under-determined system).
In both cases we removed zero rows from the system matrix.

We used Landweber's method to solve these two problems.
Figure \ref{fig:FitToNoise} shows the reconstruction errors
$\| \image^{(k)} - \gt \|_2$ and the
norms $\|\residual^{(k)}\|_2$ versus $k$, together with
the threshold $\eta\, \sqrt{m}$ and the function $\eta\,\sqrt{m-t_k}$,
i.e., here we use $\tau=1$.
The graphs confirm the monotonic decrease of the residual norm.
For both problems, the ``fit-to-noise-level'' stopping rule
terminates the iterations close to the optimal number of iterations.
A stopping rule involving $\eta\, \sqrt{m}$, on the other hand, would
terminate the iterations much too early. \hfill $\Box$

\subsection{Minimization of the Prediction Error -- UPRE}
\label{subsec:UPRE}

The key idea is to find
the number of iterations that minimizes the prediction error, i.e., the
difference between the noise-free data $\nfd = \sm\,\gt$ and the
predicted data $\sm\,\image^{(k)}$.
Statisticians refer to various measures of this difference as the
predictive risk, and the resulting method for choosing $k$ is often
called the \emph{unbiased predictive risk estimation (UPRE)} method.

Here we present the results specifically in the framework of
iterative reconstruction methods and using the matrix $\sm_k^{\#}$
defined in Eq.~\eqref{eq:Aksharp}.
Following \cite[\S 7.1]{Vogel02}, where all the details can be found,
the expected squared norm of the prediction error (the risk) is
  \begin{align*}
    \EE\bigl( \| \nfd - \sm\,\image^{(k)} \|_2^2 \bigr) =
    &\| ( \II - \sm\,\sm_k^{\#} )\, \nfd \|_2^2 \ + \\[1mm]
    &\eta^2\,\mathrm{trace}\bigl( (\sm\,\sm_k^{\#})^2 \bigr)
  \end{align*}
while the expected squared norm of the residual can be written as
  \begin{align*}
    \EE\bigl( \| \rhs \ - \ &\sm\,\image^{(k)} \|_2^2 \bigr) =
      \| ( \II - \sm\,\sm_k^{\#} )\, \nfd \|_2^2 \ + \\[1mm]
    & \eta^2\,\mathrm{trace} \bigl( (\sm\sm_k^{\#})^2 \bigr) -
      2\eta^2\,\mathrm{trace}(\sm\sm_k^{\#}) + \eta^2\, m \ .
  \end{align*}
Combining these two equations we can eliminate one of the
trace terms and arrive at the following expression for the risk:
  \begin{align*}
    \EE\bigl( \| \nfd - \sm\,\image^{(k)} \|_2^2 \bigr) = \
    &\EE\bigl( \| \rhs - \sm\,\image^{(k)} \|_2^2 \bigr) \ + \\[1mm]
    & 2\eta^2\,\mathrm{trace}(\sm\sm_k^{\#}) - \eta^2\, m \ .
  \end{align*}

Substituting the actual squared residual norm
$\|\residual^{(k)}\|_2^2 = \| \rhs - \sm\,\image^{(k)} \|_2^2$
for its expected value, we thus define the UPRE risk as a function of $k$:
  \begin{equation}
  \label{eq:UPRE}
    U^{(k)} = \|\residual^{(k)}\|_2^2 + 2\,\eta^2\, t_k - \eta^2\, m
  \end{equation}
with $t_k$ given by \eqref{eq:tk}.
A minimizer of $U^{(k)}$ will then give an approximation to a minimizer
of the prediction error.
We note that $U^{(k)}$ may not have a unique minimizer, and we therefore
choose the smallest $k$ at which $U^{(k)}$ has a local minimum.
Thus we arrive at the following stopping rule:

\medskip\noindent
\framebox[0.48\textwidth][l]{\vbox{
  \begin{tabbing}
  xxx \= xxx \= xxx \kill
  \> \underline{Stop rule:\ UPRE\,} \\[2mm]
  \> Minimize
        $U^{(k)} = \|\residual^{(k)}\|_2^2 + 2\,\eta^2\, t_k - \eta^2\, m \ .$
  \end{tabbing}
}}

\medskip

\subsection{Another Rule Based on the Prediction Error -- GCV}

The above UPRE stopping rule depends on an
estimate of the standard deviation $\eta$ of the noise -- which may or
may not be a problem in practise.
We shall now describe an alternative method for minimization of the prediction
error, derived by Wahba \cite{Wahba90},
that does not depend on knowledge of~$\eta$.

The outset for this method is the principle of cross validation.
Assume that we remove the $i$th element $b_i$ from the right-hand side
(the noisy data), compute a reconstruction $\image_{[i]}^{(k)}$, and then
use this reconstruction to compute a prediction
$\hat{b}_i = \row_i^T \image_{[i]}^{(k)}$ of the missing data~$b_i$,
where
  \[
    \row_i^T = \sm(i\, ,{:}) = \hbox{$i$th row of $\sm$.}
  \]
The goal would then be to choose the number of iterations $k$ that
minimizes the mean of all the squared prediction errors:
  \[
    \widehat{G}^{(k)} = \frac{1}{m} \, \sum_{i=1}^m \bigl( b_i - \hat{b}_i \bigr)^2 =
    \frac{1}{m} \, \sum_{i=1}^m \Bigl( b_i - \row_i^T \image_{[i]}^{(k)}
    \Bigr)^2 .
  \]
Then it is proved in \cite[Thm.~4.2.1]{Wahba90} that we can avoid the vectors
$\image_{[i]}^{(k)}$ and write $\widehat{G}^{(k)}$ directly in terms of
$\image^{(k)}$:
  \begin{equation}
    \widehat{G}^{(k)} = \frac{1}{m} \, \sum_{i=1}^m \left( \frac{b_i -
    \row_i^T \image^{(k)}}{1 - \alpha_{i}^{(k)}} \right)^2 ,
  \end{equation}
where $\alpha_i^{(k)}$ is the $i$th diagonal element of the
matrix product $\sm\,\sm_k^{\#}$ associated with~$\image^{(k)}$.

At this stage, recall that the 2-norm is invariant under an
orthogonal transformation, of which a permutation is a special case.
Specifically, if $\boldsymbol{Q}$ is an orthogonal matrix then
  \[
    \| \boldsymbol{Q} \, (\sm\,\image - \rhs) \|_2 = \| \sm\,\image - \rhs \|_2
  \]
which means that the reconstruction $\image^{(k)}$ is invariant to
such a transformation.
Unfortunately it can be proved \cite{Wahba90} that the minimizer of
$\widehat{G}^{(k)}$ is not invariant to an orthogonal transformation of the data.
In particular, it is inconvenient that a stopping rule based on
$\widehat{G}^{(k)}$ would produce a $k$ that depends on the particular
ordering of the data.

The \emph{generalized cross validation (GCV)} method circumvents this
problem by replacing all $\alpha_i^{(k)}$ with their average
  \[
    \mu^{(k)} = \frac{1}{m} \, \sum_{i=1}^m \alpha_i^{(k)} =
    \frac{1}{m}\, \mathrm{trace}(\sm\,\sm_k^{\#}) = \frac{t_k}{m} \ ,
  \]
leading to the modified measure
  \begin{align}
  \nonumber
    \widetilde{G}^{(k)} &= \frac{1}{m} \, \frac{1}{(1-\mu^{(k)})^2} \,
      \sum_{i=1}^m \bigl( b_i - \row_i^T \image^{(k)} \bigr)^2 \\
    &= \frac{\| \rhs - \sm\,\image^{(k)} \|_2^2}{m(1-t_k/m)^2} =
    m\, \frac{ \| \residual^{(k)} \|_2^2}{(m - t_k)^2} \ .
  \end{align}
The minimizer of $\widetilde{G}^{(k)}$ is, of course, independent of the
factor $m$ and hence we choose to define the GCV risk as a function of
$k$ as
  \begin{equation}
  \label{eq:GCV}
    G^{(k)} = \frac{\| \residual^{(k)} \|_2^2}{(m - t_k)^2} \ .
  \end{equation}
We have thus arrived at the following $\eta$-free stopping rule where again, in
practice, we need to estimate the quantity $t_k$:

\medskip\noindent
\framebox[0.48\textwidth][l]{\vbox{
  \begin{tabbing}
  xxx \= xxx \= xxx \kill
  \> \underline{Stop rule:\ GCV\,} \\[2mm]
  \> Minimize
        $G^{(k)} = \| \residual^{(k)} \|_2^2 \, / \, (m - t_k)^2 \ .$
  \end{tabbing}
}}

\medskip
The above presentation follows \cite[\S 4.2--3]{Wahba90}.
A different derivation of the GCV method was presented in~\cite{GoHW79};
here the coordinate system for $\mathbb{R}^m$ is rotated such that
the corresponding influence matrix becomes a circulant matrix with
identical elements along all its diagonals.
This approach leads to the same GCV risk~$G^{(k)}$ as above.

Perhaps the most important property of the GCV stopping rule
is that the value of $k$ which minimizes $G^{(k)}$ in (\ref{eq:GCV}) is
also an estimate of the value that minimizes the prediction error.
Specifically, if $k_{\mathrm{GCV}}$ minimizes the GCV risk $G^{(k)}$ and
$k_{\mathrm{PE}}$ minimizes the prediction error
$\| \nfd - \sm\,\image^{(k)} \|_2^2$,
then it is shown in \cite[\S 4.4]{Wahba90} that
  \begin{align*}
    \EE \bigl( \| \nfd - \sm\,\image^{(k_{\mathrm{GCV}})} \|_2^2 \bigr)
    \rightarrow
    \EE \bigl( \| \nfd - \sm\,&\image^{(k_{\mathrm{PE}})} \|_2^2 \bigr) \\[1mm]
    & \hbox{for} \ m \rightarrow \infty \ .
  \end{align*}

The UPRE and GCV stopping rules have the slight inconvenience that we need to
take at least one iteration too many in order to detect a minimum of $U^{(k)}$ and
$G^{(k)}$, respectively.
In practise, this is not really a problem.
For tomography problems the iteration vector $\image^{(k)}$ does not change
very much from one iteration to the next, and hence the minimum of the
error history $\| \gt - \image^{(k)} \|_2$ is usually very flat.
Hence it hardly makes any difference if we implement the UPRE and GCV
stopping rules such that we terminate the method one iteration
(or a few iterations) after the actual minimum of $U^{(k)}$ or $G^{(k)}$.

\begin{figure}
  \centering
  \includegraphics[width=0.48\textwidth]{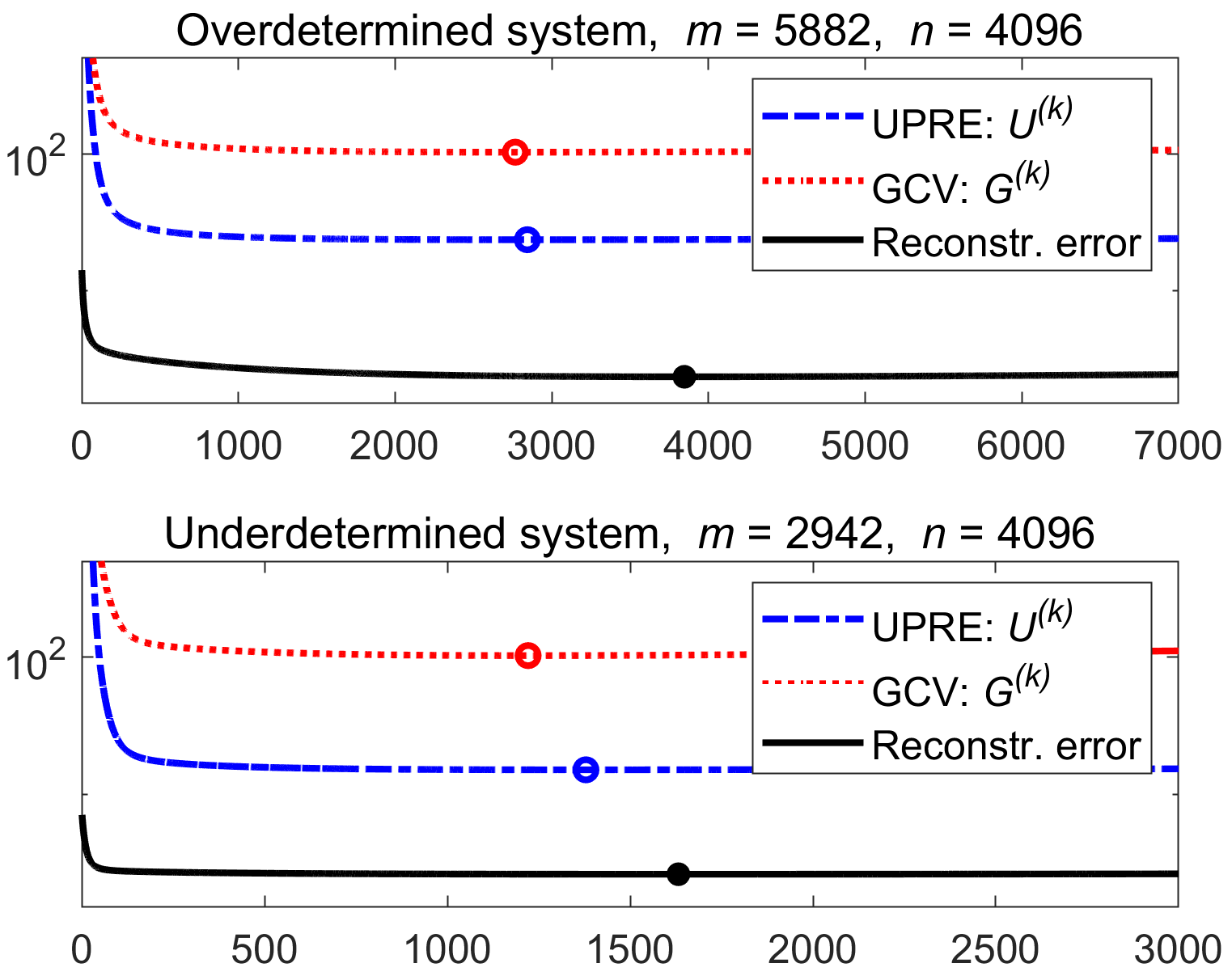}
  \caption{\label{fig:UPRE}Illustration of the UPRE and GCV stopping rules for
  Landweber's method applied to the two parallel-beam CT problems
  in Example 1 and~2.}
\end{figure}

\textbf{Example 2.}
We illustrate the UPRE and GCV stopping rules applied to Landweber's method
with the two CT problems from Example~1.
In both cases we removed zero rows from the system matrix.
Figure \ref{fig:UPRE} shows $U^{(k)}$ and $G^{(k)}$ from Eqs.\
(\ref{eq:UPRE}) and (\ref{eq:GCV})
versus $k$, together with the error histories.
The two stopping rules terminate the
iterations at approximately the same number of iterations -- not too far
from the minimum of the error history.
Note how flat the error history is:\ in practise it makes no difference if
we terminate the iterations exactly at the minimum of $U^{(k)}$ and $G^{(k)}$
or a few iterations later. \hfill $\Box$

\subsection{Stopping When All Information is Extracted --- NCP}

The above stopping rules include the trace term $t_k$ in Eq.\ (\ref{eq:tk}).
This term can be estimated at additional cost as discussed in \S\ref{subsec:trace}
below, but it is also worthwhile to consider a stopping rule that needs
neither the trace term $t_k$ nor the standard deviation $\eta$ of the noise.
The so-called NCP criterion from
\cite{HaKK06} and \cite{RuOL08} is one such method.
The considerations that underly this method are as follows:
\begin{enumerate}
\item
noisy data only contain partial information about the reconstruction,
\item
in each iteration we extract more information from the data, and
\item
eventually we have extracted all the available information in the noisy data.
\end{enumerate}
Therefore we want to monitor the properties of the residual vector.
During the initial iterations we have not yet extracted all information
present in the data and the residual still resembles a meaningful signal, while
at some stage -- when all information is extracted -- the residual starts to
appear like noise.
When we iterate beyond this point, we solely extract noise from the data
(we ``fit the noise'') and the residual vector will appear
as filtered noise where some of the noise's spectral components are removed.

To formalize this approach, in the white-noise setting of this presentation,
we need a computational approach to answering the questions:\
when does the residual vector look the most like white noise?
To answer this question, statisticians introduced the so-called
normalized cumulative periodogram.

In the terminology of signal processing,
a periodogram is identical to a discrete power spectrum defined as the
squared absolute values of the discrete Fourier coefficients.
Hence the periodogram for an arbitrary vector $\vv\in\mathbb{R}^m$ is
given by
  \begin{equation}
    \label{eq:periodogram}
    \widehat{p}_i = \bigl| \widehat{v}_i \bigr|^2 , \quad
    i = 1,2,\ldots,q \ , \qquad
    \widehat{\vv} = \mathtt{DFT}(\vv) \ .
  \end{equation}
Here, $\mathtt{DFT}$ denotes the discrete Fourier transform (computed by
means of the FFT algorithm) and $q = \lfloor m/2 \rfloor$ denotes the
largest integer such that $q \leq m/2$.
The reason for including only about half of the Fourier coefficients
in the periodogram/power spectrum is that the $\mathtt{DFT}$ of a real
vector is symmetric about its midpoint.
We then define the corresponding \emph{normalized cumulative periodogram (NCP)}
for the vector $\vv$
as the vector $\ncp(\vv)$ of length $q$ with elements, for $j = 1,2,\ldots,q$,
  \begin{equation}
    c_j(\vv) = \frac{\widehat{p}_2 + \cdots + \widehat{p}_{j+1}}
      {\widehat{p}_2 + \cdots + \widehat{p}_{q+1}}
    = \frac{\|\widehat{\vv}_{2:j{+}1} \|_2^2}{\|
      \widehat{\vv}_{2:q{+}1} \|_2^2} \ .
  \end{equation}

White noise is characterized by having a flat power spectrum (similar to white
light having equal amounts of all colors), and thus the expected value
of its power spectrum components is a constant independent of~$i$.
Consequently, the expected value of the NCP
for a white-noise vector $\vv_{\mathrm{w}}$ is the vector
  \[
    \EE\bigl(\ncp(\vv_{\mathrm{w}}\bigr)) = \boldsymbol{c}_{\mathrm{w}}
    = \left( \frac{1}{q},\frac{2}{q},\ldots,1 \right) .
  \]
How much a given vector $\vv$ deviates from being white noise can be
measured by the deviation of the corresponding $\ncp(\vv)$ from
$\boldsymbol{c}_{\mathrm{w}}$, e.g., as measured by the norm
$\| \ncp(\vv) - \boldsymbol{c}_{\mathrm{w}} \|_2$.

\begin{figure}
  \centering
  \includegraphics[width=0.48\textwidth]{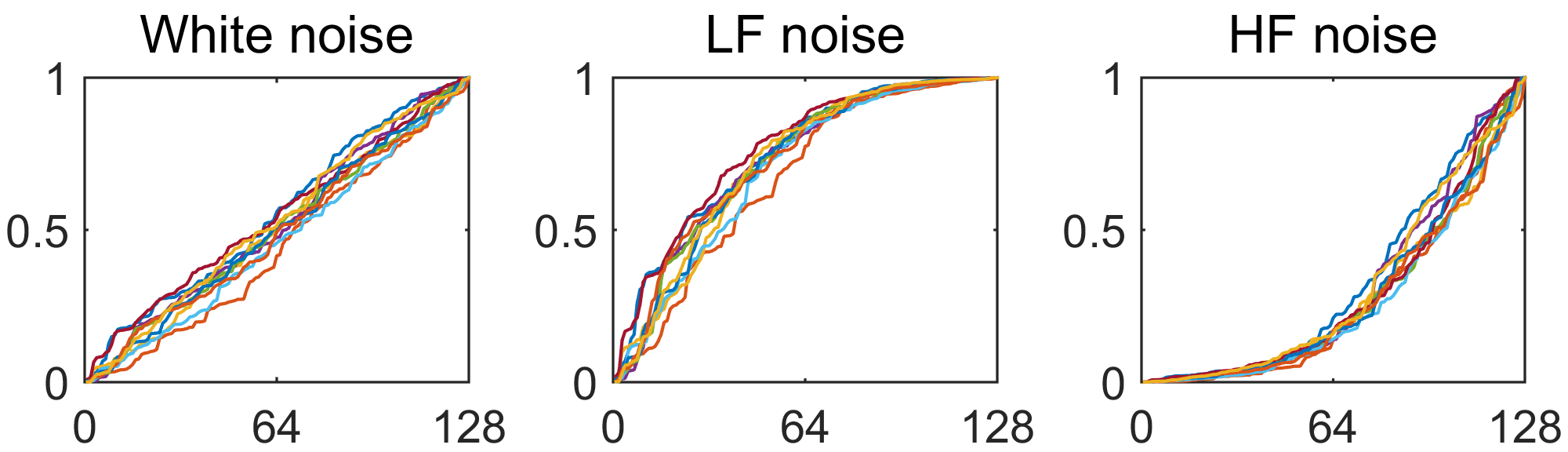}
  \caption{\label{fig:NCPs}Illustration of NCP vectors
  $\ncp(\vv)\in\mathbb{R}^{256}$ for vectors $\vv$ that are white noise (left),
  dominated by low-frequency components (middle), and dominated by
  high-frequency components (right).}
\end{figure}

\textbf{Example 3.}
Figure \ref{fig:NCPs} illustrates the appearance of NCP vectors
$\ncp(\vv)$ for vectors $\vv$ of length $m=256$
with different spectra.
The completely flat spectrum for white noise corresponds to a straight line
from $(0,0)$ to $(q,1)$ with $q = \lfloor 256/2 \rfloor = 128$.
The left plot shows NCPs for 10 random realizations of white noise, and they
are all close to the ideal white-noise NCP $\boldsymbol{c}_{\mathrm{w}}$.
The middle and right plots show NCPs for random vectors that are dominated by
low-frequency and high-frequency components, respectively; their systematic
deviation from $\boldsymbol{c}_{\mathrm{w}}$ is obvious. \hfill $\Box$

To utilize the NCP framework in the algebraic iterative methods for tomographic
reconstruction, a first idea might be to terminate the iterations when
the deviation measured by
$\| \ncp(\residual^{(k)}) - \boldsymbol{c}_{\mathrm{w}} \|_2$
exhibits a minimum.
However, this would be a bit naive since the residual vector does not really
correspond to a 1D signal of length $m$.
Rather, the right-hand side $\rhs$ consists of a number of projections,
one for each angle of the measurements -- and the residual vector
inherits this structure.
Hence, a better approach is to apply an NCP analysis to each
projection's residual, and then combine this information into a simple measure.

Depending on the CT scanner, each projection is either a 1D or 2D image,
when we perform 2D and 3D reconstructions, respectively.
To simplify our presentation, we assume that our data consists of
$m_{\theta}$ 1D projections, one for each projection angle
$\theta_1,\theta_2,\ldots,\theta_{m_{\theta}}$.
We also assume that the data are organized such that we can partition the
right-hand side $\rhs$ and the residual vector into $m_{\theta}$ sub-vectors,
  \begin{equation}
    \rhs = \begin{pmatrix} \rhs_1 \\ \rhs_2 \\ \vdots \\ \rhs_{m_{\theta}}
    \end{pmatrix} , \qquad
    \residual^{(k)} = \begin{pmatrix} \residual_1^{(k)} \\ \residual_2^{(k)}
    \\ \vdots \\ \residual_{m_{\theta}}^{(k)} \end{pmatrix} ,
  \end{equation}
with each sub-vector corresponding to a single 1D projection.
Now define the corresponding quantities
  \begin{equation}
    \nu_{\ell}^{(k)} = \bigl\| \ncp\bigl(\residual_{\ell}^{(k)}\bigr) -
    \boldsymbol{c}_{\mathrm{w}} \bigr\|_2 , \quad \ell = 1,2,\ldots,m_{\theta}
  \end{equation}
that measure the deviation of each residual sub-vector from being white noise.
Then for the $k$th iteration
we propose to measure the residual's deviation from being white
noise by averaging the above quantities, i.e., by means of the ``NCP-number''
  \begin{equation}
  \label{eq:NCPk}
    N^{(k)} = \frac{1}{m_{\theta}} \, \sum_{\ell=1}^{m_{\theta}} \nu_{\ell}^{(k)} \ .
  \end{equation}
This multi-1D approach for 2D reconstruction problems
leads to the following stopping rule:

\medskip\noindent
\framebox[0.48\textwidth][l]{\vbox{
  \begin{tabbing}
  xxx \= xxx \= xxx \kill
  \> \underline{Stop rule:\ NCP\,} \\[2mm]
  \> Minimize
        $\displaystyle N^{(k)} = \frac{1}{m_{\theta}} \, \sum_{\ell=1}^{m_{\theta}}
        \bigl\| \ncp\bigl(\residual_{\ell}^{(k)}\bigr) -
    \boldsymbol{c}_{\mathrm{w}} \bigr\|_2 \ .$
  \end{tabbing}
}}

\medskip
In the case of 3D reconstructions, where the data consist of a collection
of 2D images, the computation of $\nu_{\ell}^{(k)}$ should take this
into consideration.
In particular we need to define the NCP vector
$\ncp\bigl(\residual_{\ell}^{(k)}\bigr)$ when the residual
sub-vector $\residual_{\ell}^{(k)}$ represents an image;
how to do this is explained in~\cite{HaKK06}.

Similar to the previous stopping rules, in practise it is more convenient
to implement the NCP stopping rule such that we terminate the iterations
at the first iteration $k$ for which $N^{(k)}$ increases.
There is no theory to guarantee that $N^{(k)}$ will behave smoothly, and
we occasionally see that $N^{(k)}$ exhibits a minor zig-zag behavior.
Hence it may be necessary to apply the NCP stopping rule to a smoothed
version of the NCP-numbers, obtained by applying a ``local'' low-pass
filter to the $N^{(k)}$-sequence.

\begin{figure*}
  \centering
  \includegraphics[width=0.9\textwidth]{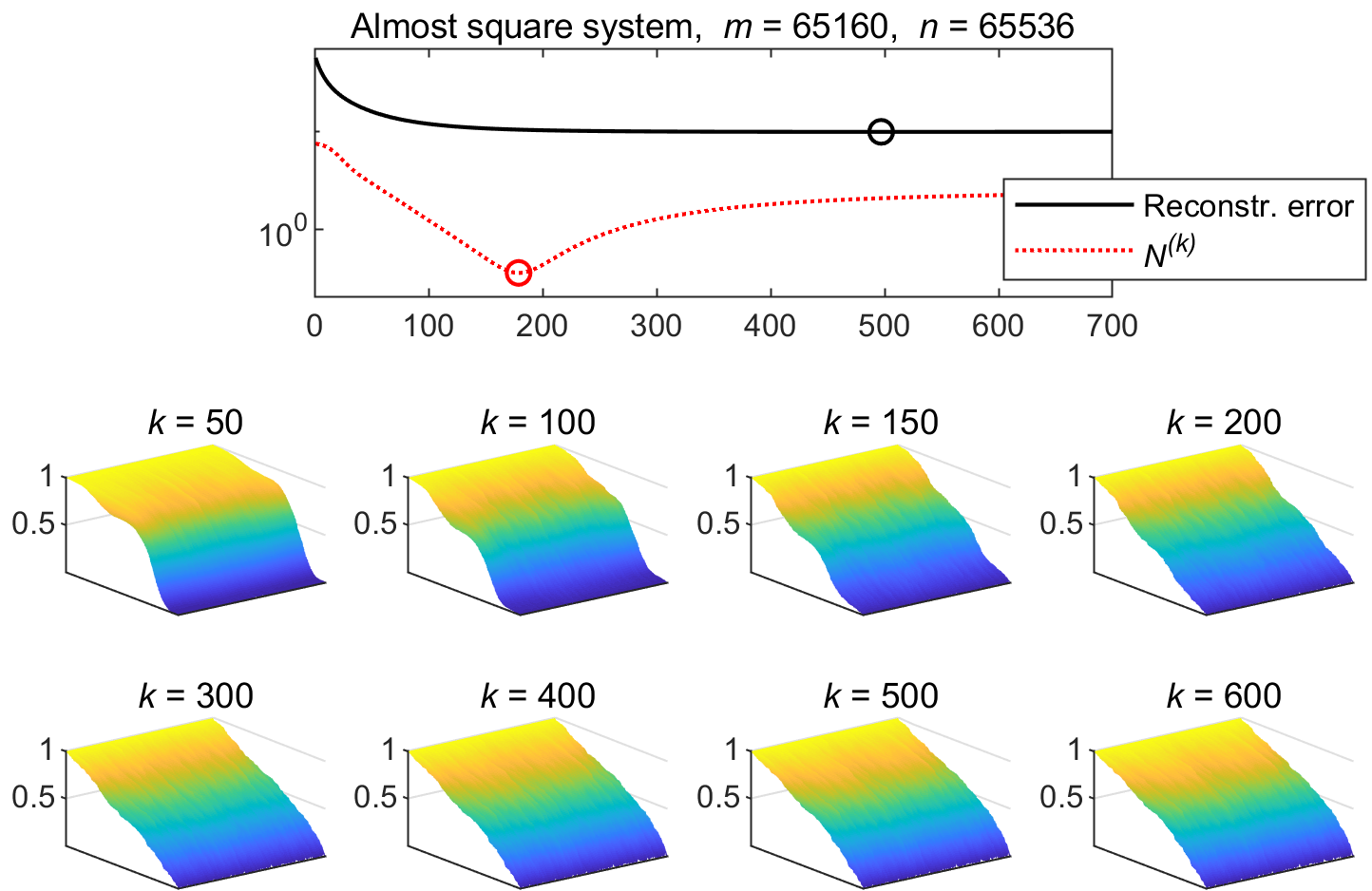}
  \caption{\label{fig:NCPagain}Illustration of the NCP stopping rules for
  Landweber's method applied to a parallel-beam test problem.
  We also show surface plots of the matrix
  $\bigl[ \ncp\bigl(\residual_1^{(k)}\bigr),
    \ncp\bigl(\residual_2^{(k)}\bigr),\ldots,
   \ncp\bigl(\residual_{m_{\theta}}^{(k)}\bigr) \bigr]$
  for selected iterations~$k$.
  This stopping rule leads to a somewhat premature termination of
  the iterations at $k_{\mathrm{NCP}}=179$
  (the minimum error occurs for $k=497$ iterations), but it
  should be noted that the error does not change much between
  iterations $179$ and $700$.}
\end{figure*}

\textbf{Example 4.}
We illustrate the NCP stopping rule with a parallel-beam CT problem
with image size $256 \times 256$ and with $362$ detector pixels and
projection angles $1^{\circ}, 2^{\circ},\ldots,180^{\circ}$.
The performance is shown in Fig.~\ref{fig:NCPagain} together with surface
plots of the matrix
$\bigl[ \ncp\bigl(\residual_1^{(k)}\bigr),
 \ncp\bigl(\residual_2^{(k)}\bigr),\ldots,
 \ncp\bigl(\residual_{m_{\theta}}^{(k)}\bigr) \bigr]$ for
selected iterations~$k$.
We clearly see the changing shape of the NCP vectors
$\ncp\bigl(\residual_{\ell}^{(k)}\bigr)$ as $k$ increases.
The minimum of $N^{(k)}$ is obtained at $k_{\mathrm{NCP}}=179$.
This is somewhat early, considering that the minimum reconstruction
error is obtained at $k=497$ iterations -- but on the other hand,
the reconstruction and the error changes only little between
iterations $150$ and $700$. \hfill $\Box$

\section{Estimation of the Trace Term} 
\label{subsec:trace}

The FTNL, UPRE and GCV stopping rules include the
term $t_k=\mathrm{trace}(\sm\,\sm_k^{\#})$.
To make these methods practical to use, we need to be able to
estimate this trace term efficiently, without having to
compute the SVD of the system matrix $\sm$ or
form the influence matrix~$\sm\,\sm_k^{\#}$.
The most common way to compute this estimate is via a Monte Carlo approach.

Underlying this approach is the following result from~\cite{Girard89}.
If $\overline{\ww}\in\mathbb{R}^m$ is a random vector with elements
$\overline{w}_i \sim \mathcal{N}(0,1)$,
and if $\boldsymbol{S} \in \mathbb{R}^{m \times m}$ is a symmetric matrix,
then $\overline{\ww}^T \boldsymbol{S}\, \overline{\ww}$ is an unbiased estimate
of~$\mathrm{trace}(\boldsymbol{S})$.
Therefore
$\overline{t}_k^{\,\mathrm{est}} = \overline{\ww}^T \sm\,\sm_k^{\#}\, \overline{\ww}$
is an unbiased estimator of $t_k = \mathrm{trace}(\sm\,\sm_k^{\#})$.

To compute this estimate we need to compute the matrix-vector
product $\sm_k^{\#} \overline{\ww}$ efficiently.
Recalling the definition of $\sm_k^{\#}$ in Eq.~(\ref{eq:Aksharp}),
this can be done simply by applying the algebraic iterative method
to the system $\sm\,\overline{\boldsymbol{\xi}} = \overline{\ww}$
which, after $k$ iterations, produces the iteration vector
$\overline{\boldsymbol{\xi}}^{(k)} = \sm_k^{\#} \overline{\ww}$.
The resulting estimate
  \begin{equation}
  \label{eq:tkestoriginal}
    \bar{t}_k^{\,\mathrm{est}} = \overline{\ww}^T \sm\,\,
    \overline{\boldsymbol{\xi}}^{(k)} = ( \sm^T \overline{\ww})^T \,
    \overline{\boldsymbol{\xi}}^{(k)}
  \end{equation}
is the standard Monte Carlo trace estimate from~\cite{Girard89}.
In an efficient implementation of (\ref{eq:tkestoriginal}) the
vector $\sm^T \overline{\ww}$ is pre-computed and stored.

An alternative approach was presented in \cite{SantosDePierro03}.
This approach also applies to the general method in \eqref{eq:SIRT}
with $\DD=\II$ and with a general $m\times m$ matrix $\MM$
(it is not required to be symmetric).
When we apply such a method with an arbitrary nonzero starting vector
$\boldsymbol{\xi}^{(0)}$ to the system
$\sm\,\boldsymbol{\xi} = \boldsymbol{0}$, then it follows from
Eq.~(\ref{eq:Aksharp}) that the iterates are
  \[
    \boldsymbol{\xi}^{(k)} =
    (\II - \relaxp\,\sm^T \boldsymbol{B}\,\sm)^k\,\boldsymbol{\xi}^{(0)} \ .
  \]

Then it is shown in \cite{SantosDePierro03} that if we use a random starting vector
$\boldsymbol{\xi}^{(0)} = \ww\in\mathbb{R}^n$
with elements $w_i \sim \mathcal{N}(0,1)$,
and if $\boldsymbol{\xi}^{(k)}$ denotes the corresponding iterations
for the system $\sm\,\boldsymbol{\xi} = \boldsymbol{0}$, then
$\ww^T \boldsymbol{\xi}^{(k)}$
is an unbiased estimator of $n-\mathrm{trace}(\sm\,\sm_k^{\#})$.
This leads to the alternative trace estimate
  \begin{equation}
  \label{eq:tkest}
    t_k^{\,\mathrm{est}} = n - \ww^T \boldsymbol{\xi}^{(k)} .
  \end{equation}

In order to use either of these trace estimates instead of the exact $t_k$,
we must simultaneously apply the iterative method to two right-hand sides,
which essentially doubles the amount of work.
The Landweber method with the two different trace estimation schemes
are shown below.

If we are willing to increase the overhead further, we can
compute a more robust estimate of $t_k$ by applying the above idea
to several random vectors and computing the mean or median
of the $t_k^{\,\mathrm{est}}$-values.

\medskip\noindent
\framebox[0.48\textwidth][l]{\vbox{
  \begin{tabbing}
  xxx \= xxx \= xxx \kill
  \> \underline{Landweber method with \eqref{eq:tkestoriginal} trace estimator}
     \\[2mm]
  \> $\overline{\ww}$ = random $m$-vector for trace estimation \\
  \> $\image^{(0)}$ = initial vector \\
  \> $\overline{\boldsymbol{\xi}}^{(0)} = \boldzero$ initial zero vector
     for trace estimation \\
  \> $\vz = \sm^T \overline{\ww}$ \\
  \> for $k=0,1,2,\ldots$ \\[1mm]
  \> \> $\image^{(k+1)} = \image^{(k)} + \relaxp\,\sm^T (\rhs - \sm\,\image^{(k)} )$
     \\[1mm]
  \> \> $\overline{\boldsymbol{\xi}}^{(k+1)} = \overline{\boldsymbol{\xi}}^{(k)}
        + \relaxp\,\sm^T (\overline{\ww} - \sm\,\overline{\boldsymbol{\xi}}^{(k)} )$
     \\[1mm]
  \> \> $\bar{t}_{k+1}^{\,\mathrm{est}} = \vz^T \overline{\boldsymbol{\xi}}^{(k+1)}$
        trace estimate \\[1mm]
  \> \> stopping rule goes here \\
  \> end
  \end{tabbing}
}}

\medskip\noindent
\framebox[0.48\textwidth][l]{\vbox{
  \begin{tabbing}
  xxx \= xxx \= xxx \kill
  \> \underline{Landweber method with \eqref{eq:tkest} trace estimator}
     \\[2mm]
  \> $\ww$ = random $n$-vector \\
  \> $\image^{(0)}$ = initial vector \\
  \> $\boldsymbol{\xi}^{(0)} = \ww$ initial vector for
     for trace estimation \\
  \> for $k=0,1,2,\ldots$ \\[1mm]
  \> \> $\image^{(k+1)} = \image^{(k)} + \relaxp\,\sm^T (\rhs - \sm\,\image^{(k)} )$
     \\[1mm]
  \> \> $\boldsymbol{\xi}^{(k+1)} = \boldsymbol{\xi}^{(k)}
        + \relaxp\,\sm^T (\boldzero - \sm\,\boldsymbol{\xi}^{(k)} )$
     \\[1mm]
  \> \> $t_{k+1}^{\,\mathrm{est}} = n - \ww^T \boldsymbol{\xi}^{(k+1)}$
        trace estimate \\[1mm]
  \> \> stopping rule goes here \\
  \> end
  \end{tabbing}
}}

\medskip

\begin{figure}
  \centering
  \includegraphics[width=0.48\textwidth]{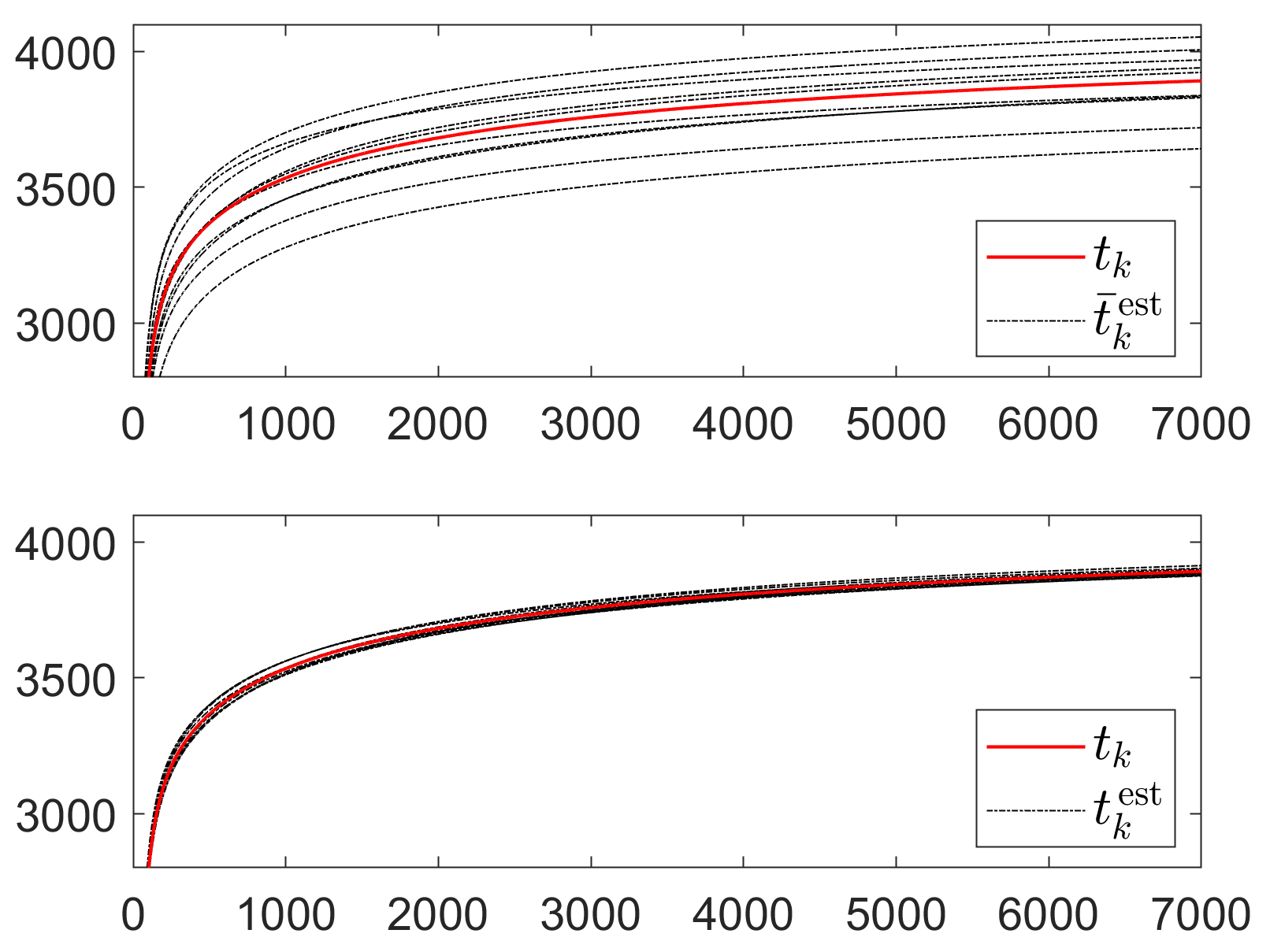}
  \caption{\label{fig:TraceEstOnly}Comparison of the two trace estimates
  $\bar{t}_k^{\,\mathrm{est}}$ and $t_k^{\,\mathrm{est}}$ for
  Landweber's method applied to the over-determined test problem
  from Example~1.
  The thick red line is the exact trace $t_k$, and the thin black lines
  are the trace estimates for 10 different random vectors
  $\overline{\ww}$ and $\ww$.}
\end{figure}

\textbf{Example 5.}
We illustrate the two trace estimates $\bar{t}_k^{\,\mathrm{est}}$ and
$t_k^{\,\mathrm{est}}$ for Landweber's method applied to the over-determined
test problem from Example~1.
Figure \ref{fig:TraceEstOnly} shows the trace estimates for 10
different realizations of the random vectors $\overline{\ww}$ and $\ww$,
together with the exact trace~$t_k$.
We see that the estimate $t_k^{\,\mathrm{est}}$, shown in the bottom plot,
has the smallest variance.
We are not aware of theory that supports this observation. \hfill $\Box$

\begin{figure}
  \centering
  \includegraphics[width=0.48\textwidth]{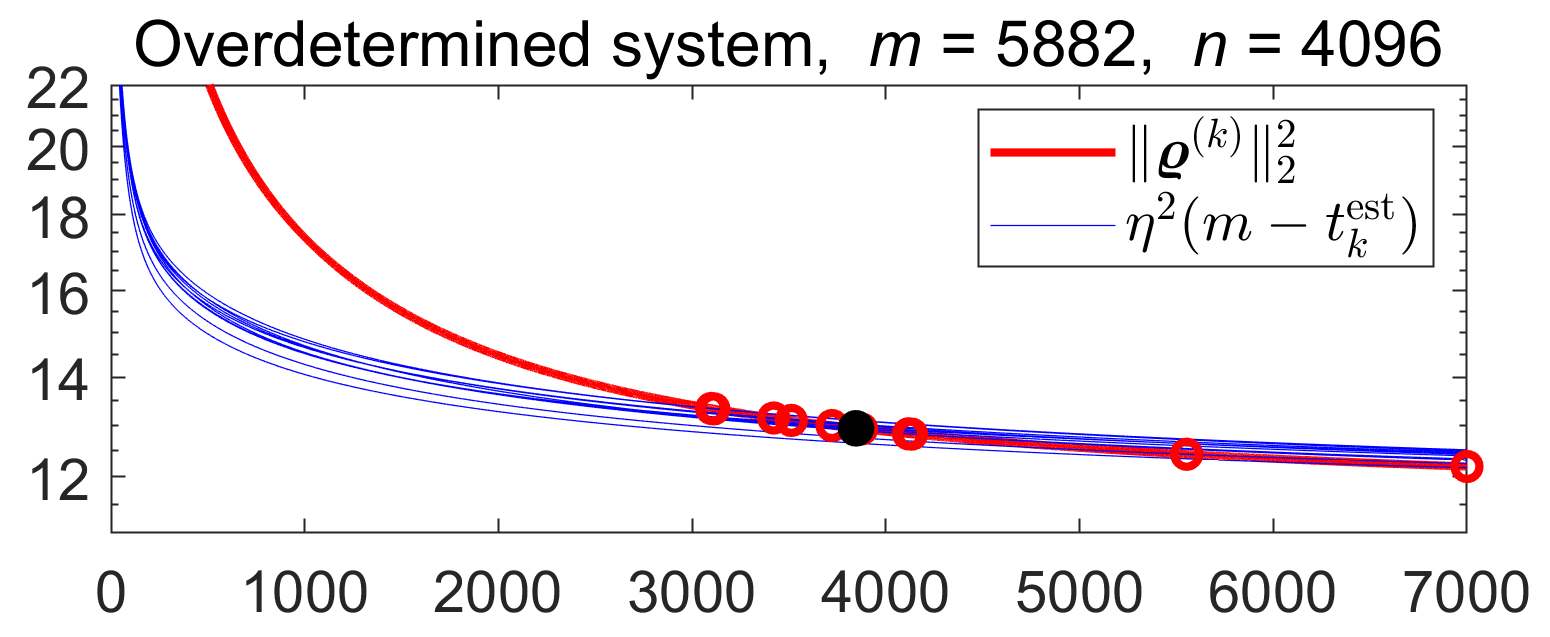}
  \caption{\label{fig:TraceEstFitToNoise}Illustration of the use
  of the trace estimate $t_k^{\,\mathrm{est}}$ in the FTNL stopping rule
  for Landweber's method applied to the over-determined test problem
  from Examples 1.
  We used 10 different random vectors $\ww$ in (\ref{eq:tkest}) and
  the corresponding 10 intersections between $\|\residual^{(k)}\|_2^2$
  (thick red line) and $\eta^2\,(m-t_k^{\,\mathrm{est}})$ (thin blue lines)
  are shown by the red circles.
  The black dot shows the intersection with the exact $\eta^2\,(m-t_k)$.}
\end{figure}

\textbf{Example 6.}
Continuing from the previous example, Fig.~\ref{fig:TraceEstFitToNoise}
illustrates the use of the trace estimate $t_k^{\,\mathrm{est}}$ in the
FTNL stopping rule.
To show the variability of the stopping rule we used 10 different random
vectors $\ww$,
leading to 10 different realizations of $\eta^2\,(m-t_k^{\,\mathrm{est}})$.
Their intersections with $\|\residual^{(k)}\|_2^2$ are shown by the
red circles, corresponding to stopping the iterations at
  \begin{align*}
    k = 3100,\ &3112,\ 3421,\ 3512,\ 3722,\\ &3875,\ 4117,\ 4133,\ 5553,\ 7000 .
  \end{align*}
The black dot marks the intersection of the exact of $\|\residual^{(k)}\|_2^2$
with $\eta^2\,(m-t_k)$, corresponding to iteration $k=3846$. \hfill $\Box$

\section{Large-Scale Numerical Example}
\label{sec:large}

In this section we use a simulated large-scale CT reconstruction problem
to illustrate the use of the GCV and NCP stopping rules described above.
We focus on an application in dynamic tomography where
the time scale of the process being examined dictates the use of
a small number of projections as well as short exposure times of each projection.
This leads to challenging reconstruction problems where it is critical to
use a stopping rules that terminates the iterations such that
$\image^{(k)}$ is as close as possible to $\bar{\bm{x}}$
and without having knowledge of the noise level in the data.

Specifically we study the use of
the GCV and NCP stopping rules applied to the reconstruction of
a single time step in a simulation of a dynamic CT experiment.
The dynamic process under study is the
separation of an emulsion of oil and water in a porous rock;
the components separate vigorously
over time, due to the two fluids being immiscible.

The basis of our simulation is a segmentation of a nano-CT scan of a piece of
chalk from the Hod field in the North Sea Basin
(sample id HC \#15) which was scanned, reconstructed and segmented as
described in \cite{AWR,APL}.
A subset consisting of $200\times256\times256$ voxels is chosen for the fluid simulation.
Pixels outside a radius of 124 pixels from the center axis are set to zero
to form a cylinder, which is mirrored along its vertical axis to ensure that the
multiphase flow simulation has periodic boundary conditions.
The flow simulation is done with a phase-field Lattice Boltzmann method for systems that
are isothermal and incompressible \cite{Fakhari2017a, Fakhari2017b}.
The simulation produces phase values for each voxel that describes the
fraction of oil and water in the voxel.
These phase values are converted
to attenuation values based on values measured in \cite{Lin2018} where a
sandstone filled with a brine and oil is imaged with X-rays at 80~keV.

\begin{figure}
\centering
\includegraphics[keepaspectratio, height=0.25\paperheight]{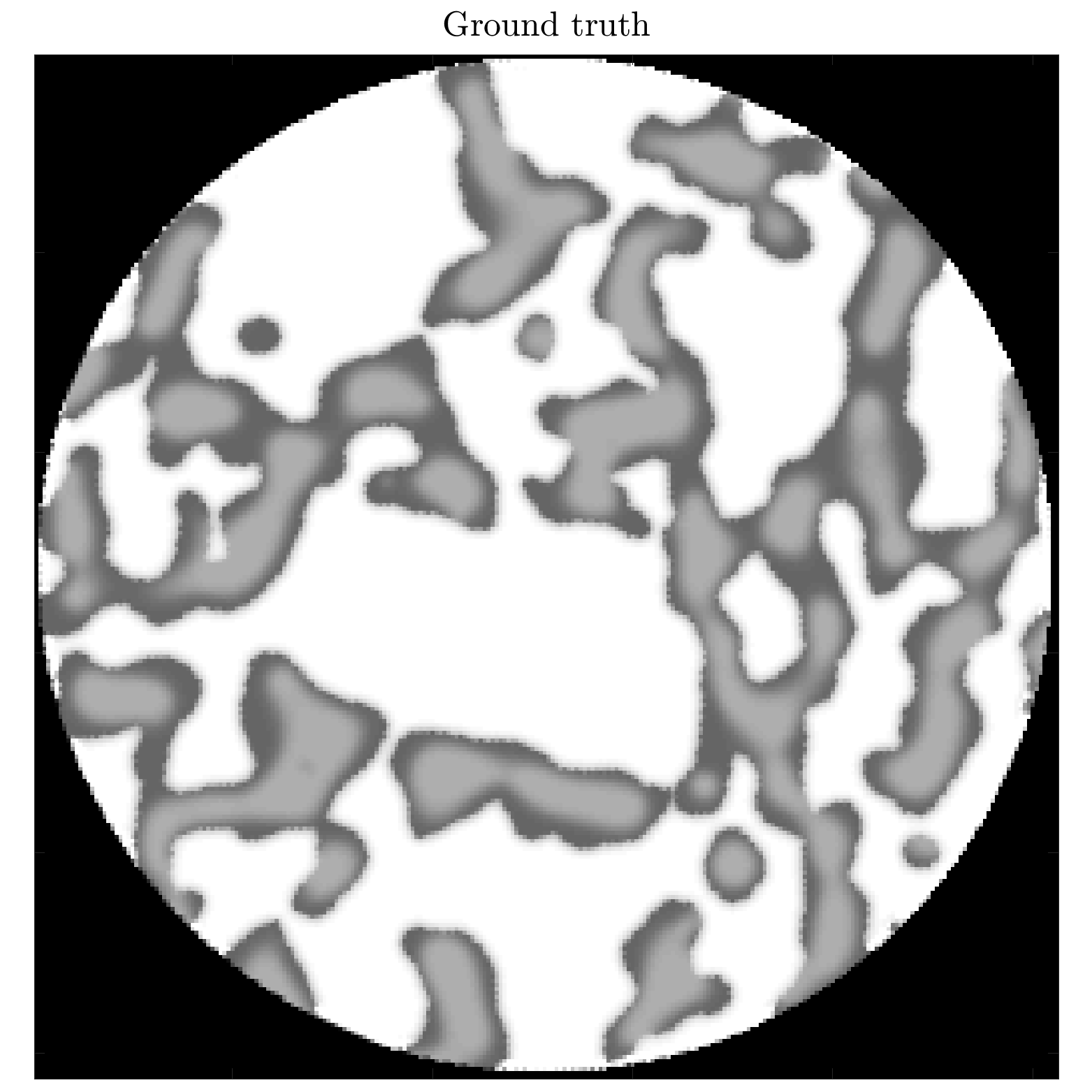}
\caption{A single slice of the volume being examined. White corresponds to
rock while light and dark grey correspond to water and oil respectively.}
\label{fig:ground_truth}
\end{figure}

Figure \ref{fig:ground_truth} shows
a slice from a single time frame in the simulation;
the rock matrix is white, the water phase is light grey, and the
oil phase is dark grey.
The time frame which was chosen for testing the stopping rules
is fairly early in the simulation where multiple interfaces between the two
fluids are present.

\subsubsection*{Forward projection}
The forward projection of the volume is performed using the ASTRA toolbox
with a parallel beam geometry \cite{ASTRA2011,ASTRA2015,ASTRA2016}.
We use 362 detector pixels and 360 projection angles.
The forward projection in ASTRA can be considered an ideal experiment with
monochromatic X-rays and infinite brilliance, i.e., without any noise.

\subsubsection*{Noise}
We create noisy data from the above clean data in such a way that we
emulate the noise present in
X-ray tomography as a result of the finite count of photons,
cf.\ \cite[\S 4.4]{SIAMbook}.
Specifically, if $\nfd = \sm\,\gt$ denotes the clean data computed
by means of ASTRA, then the corresponding X-ray intensities at the detector
are given by
  \[
    \bar{I}_i = I_0 \exp(-\bar{b}_i) \ , \quad i=1,\ldots,m \ ,
  \]
where $I_0$ is the source's intensity.
We then use $\bar{I}_i$ as the ex\-pected value in a Poisson distribution to obtain
noisy intensities
  \[
    I_i = \mathcal{P}(\bar{I}_i) = \mathcal{P}(I_0 \exp(-\bar{b}_i)) \ , \quad i=1,\ldots,m \ .
  \]
Finally, we convert these noisy intensities back to the noisy data vector
$\rhs$ via the relation
  \[
    b_i = -\log(I_i/I_0) \ , \quad i=1,\ldots,m \ .
  \]

We use three different noise levels $0.25\%$, $1\%$
and $5\%$ which visually corresponds to low, moderate and high noise.
The noise level is given by
  \begin{align}
    \rho = \frac{\| \noise \|_2}{\| \nfd \|_2} \ , \qquad
    \noise = \rhs - \nfd \ ,
  \end{align}
where $\bm{e}$ denotes the measurement error in Eq.~\eqref{eq:Ax}.
This noise does not exactly fit with the assumption of white Gaussian noise
which is used for the previous derivations, but it is a good approximation
to the noise present in CT experiments.

\subsubsection*{Reconstruction}
We compute reconstructions from the simulated
projection data with the ASTRA toolbox by using the
Simultaneous Iterative Reconstruction Technique (SIRT) iterative method.
This is a special case of the general method in Eq.~\eqref{eq:SIRT}
where the diagonal matrices $\DD$ and $\MM$ contain the inverse column and row sums of
$\bm{A}$:
  \[
    d_{jj}=1/\sum_{i=1}^m a_{ij} \quad \hbox{and} \quad
    m_{ii}=1/\sum_{j=1}^n a_{ij} \ .
  \]

We perform 1000 SIRT iterations
and for UPRE and GCV the trace $t^{(k)}$ in $U^{(k)}$ and $G^{(k)}$
is estimated using Eq.~\eqref{eq:tkestoriginal}.
Note that only a single random vector $\bm{w}$ is used to reduce computation time.
This is not a concern in this specific case as the estimation of $t^{(k)}$ proved
very stabile for different random vectors.

The forward projection of $\image^{(k)}$ used in the calculation of $N^{(k)}$
is computed with ASTRA, and the remaining part of the algorithm is calculated
with CuPy, a Python package which makes it possible to offload calculations
to a CUDA-compatible GPU to improve the computation time of $N^{(k)}$.
The vector $\vv$ in \eqref{eq:periodogram} is
padded with zeros such that its length can be written in the form
$n = 2^a + 3^b + 5^c + 7^d$ as this substantially speeds up the calculation
when the $\mathtt{DFT}$ is calculated with CUDA.

\subsubsection*{Results}
As previously mentioned, we reconstruct the data at three different noise levels.
Moreover, we subsample the number of projections used for the
reconstruction such that it is performed with 360, 120 and 45 projection angles.
This leads to 9 different data sets.

\begin{figure*}
\centering
    \includegraphics[keepaspectratio, width=\textwidth]{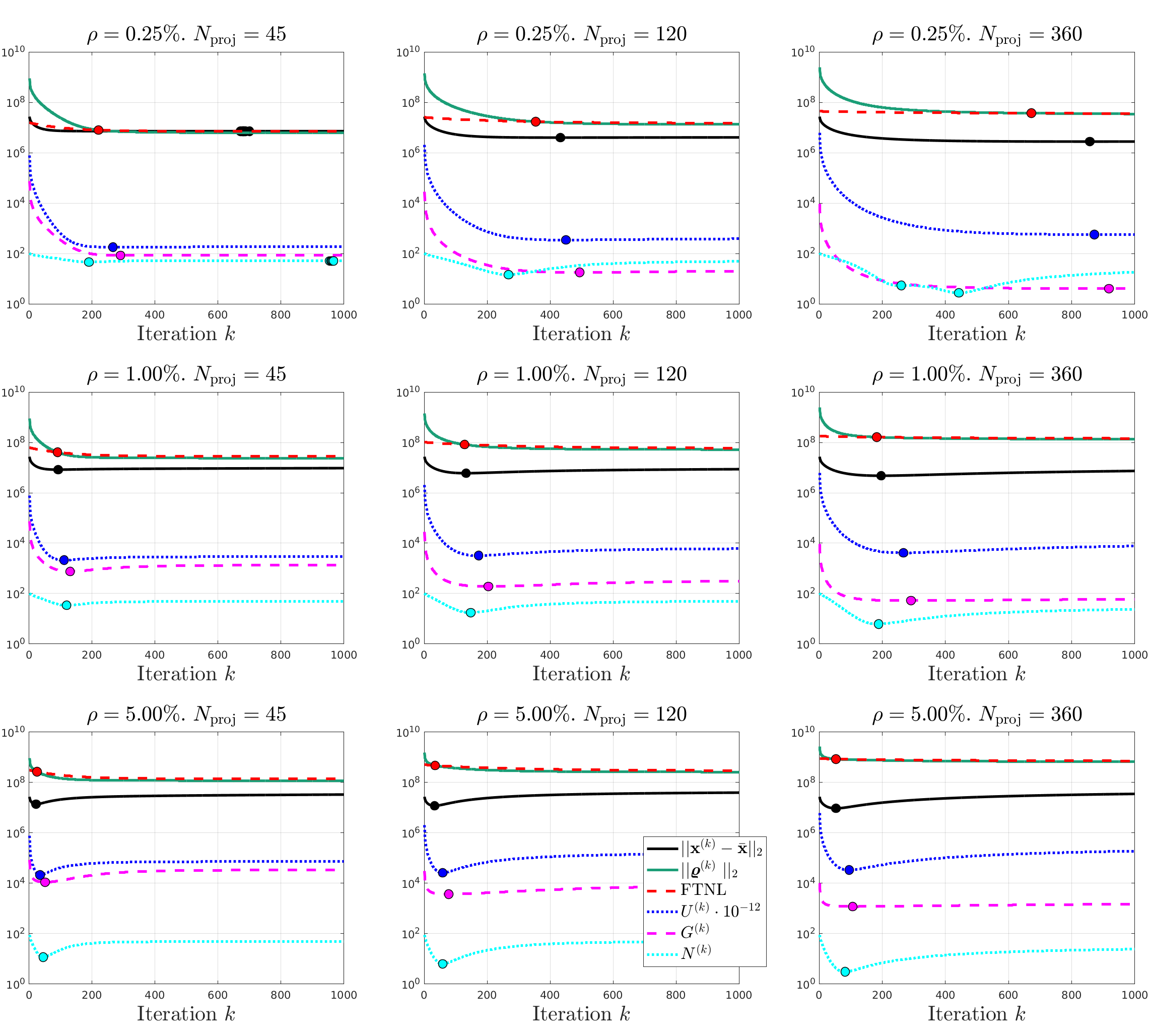}
    \caption{
    Illustration of the four different stopping rules for the large-scale example.
    The filled circles on each curve represent the minimum.
    Each column has a varying number of projections and each row
    has a varying amount of noise, as shown in the titles of each subplot.}
    \label{fig:stopping_rules}
\end{figure*}

Figure \ref{fig:stopping_rules} shows graphs of
$\|\residual^{(k)} \|_2$, $\tau \, \eta\,\sqrt{m-t_k}$ (called FTNL),
$U^{(k)}$, $G^{(k)}$, $N^{(k)}$ and $\| \image^{(k)}-\gt \|_2$
along with their local minima.
In general we see that the FTNL and UPRE stopping rules perform well in this
simulated example; but they depend on knowledge of the noise level $\eta$
which is rarely known for real data.
An advantage of the GCV and NCP stopping rules is that they do not rely
on an estimate of~$\eta$.
GCV performs well in the case with 1\% and 5\% noise,
but it overestimates the number of iterations with the small noise level $\rho=0.25\%$.
The NCP stopping rule is also a robust method and it performs well
for all noise levels.
It is worth noting that the reconstruction error is very flat for this noise level,
which means the exact amount of iterations used is less critical.

  \begin{figure*}
    \centering
    \includegraphics[keepaspectratio, width=\textwidth]{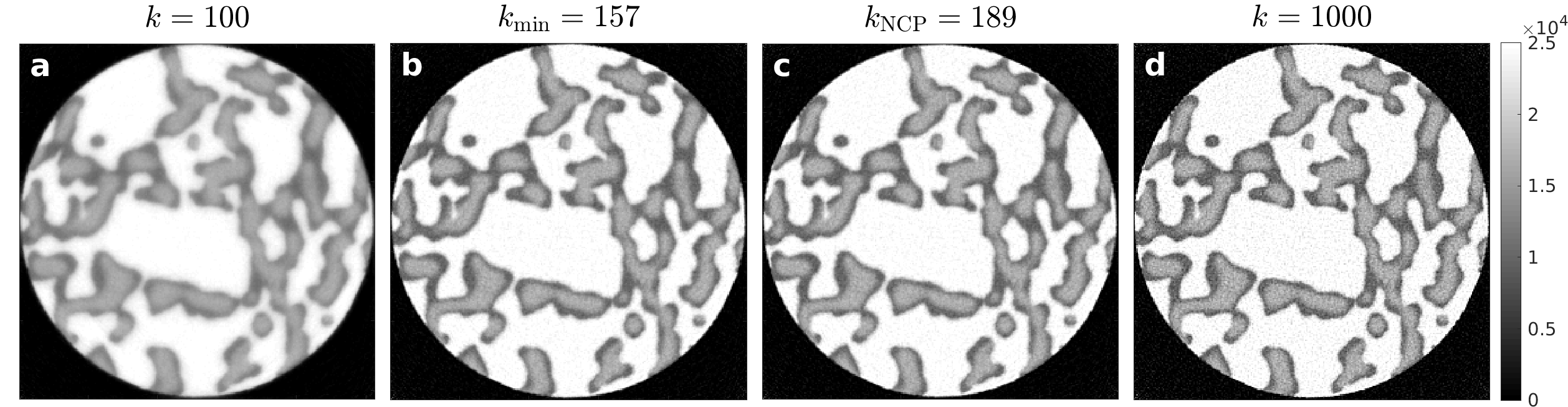}
    \caption{Illustration of semi-convergence for the large-scale example with
    $360$ projections and noise level $\rho=1\%$.
    Results for four different iteration numbers are shows.
    Image \textbf{a} is the reconstruction after
    100 iterations. Image \textbf{b} is the reconstruction at $k_{\min} = 157$
    which is the number of iterations that minimizes the reconstruction error.
    Image \textbf{c} is the reconstruction at $k_{\mathrm{NCP}} = 189$ which is
    the number of iterations that minimizes $N^{(k)}$. Finally, image
    \textbf{d} is the reconstruction after $k=1000$ iterations.}
    \label{fig:optimal_vs_actual}
  \end{figure*}

Figure \ref{fig:optimal_vs_actual} shows the effect of semi-convergence on the
data set when it is reconstructed with $N_{\mathrm{proj}}=360$ and
$\rho=1.00\%$. A single slice of the reconstruction is shown in all subfigures.
All images in Fig.~\ref{fig:optimal_vs_actual} are truncated such that their
intensities are between 0 and 25,000.
Image \textbf{a} is the reconstruction after $k=100$ iteration where
it still has a blurred appearance, showing that more iterations are necessary.
Images \textbf{b} and \textbf{c} are
the reconstructions at the number of iterations which minimize the
reconstruction error and $N^{(k)}$, respectively.
The appearance of these reconstructions is very
similar, but \textbf{c} has a slight increase in noise.
The rightmost image
\textbf{d} is the reconstruction after 1000 iterations where it is
noticeable more noisy than \textbf{b} and \textbf{c}.

\section{Conclusion}

We surveyed several state-of-the-art stopping rules, based on statistical
considerations, that are useful
for algebraic iterative reconstruction methods in X-ray computed tomography (CT)\@.
Common for these stopping rules is that they seek to terminate the iterations
at the optimal point where the reconstruction error and the noise error
balance each other.
They are easy to use and they are also easy to integrate in existing software.
We also illustrated the use of two of these methods for a large-scale CT problem
related to the study of multiphase flow in chalk.
Our numerical experiments show that especially the NCP stopping rule -- which is
based on statistical properties of the residual and does not depend on knowledge
of the noise level -- works well for this problem.

\section*{Acknowledgement}
We thank the referees for valuable comments.


\end{document}